\documentclass[12pt]{article}
\usepackage{amsmath,amsfonts,amssymb,amscd}
\usepackage{graphicx}
\usepackage[T2A]{fontenc}

\newtheorem{lem}{Lemma}[section]
\newtheorem{prop}{Proposition}[section]

\newtheorem{cor}{Corollary}[section]
\newtheorem{rem}{Remark}[section]

\newtheorem{dfn}{Definition}[section]

\makeatletter \@addtoreset{equation}{section} \makeatother

\newcommand{\mC}{\mathbb{C}}

\newcommand{\mR}{\mathbb{R}}
\newcommand{\mT}{\mathbb{T}}
\newcommand{\mZ}{\mathbb{Z}}
\newcommand{\mN}{\mathbb{N}}

\newcommand{\one}{{\bf 1}}

\newcommand{\bX}{{\bf X}}

\newcommand{\bc}{{\bf c}}

\newcommand{\bg}{{\bf g}}

\newcommand{\bT}{{\bf T}}

\newcommand{\calB}{{\cal B}}

\newcommand{\calH}{{\cal H}}

\newcommand{\calM}{{\cal M}}

\newcommand{\calO}{{\cal O}}

\newcommand{\calR}{{\cal R}}

\newcommand{\calW}{{\cal W}}
\newcommand{\calX}{{\cal X}}

\newcommand{\AACF}{{\cal AACF}}

\newcommand{\ACF}{{\cal ACF}}
\newcommand{\DT}{{\cal DT}}

\newcommand{\frakh}{{\mathfrak h}}
\newcommand{\frakX}{{\mathfrak X}}

\newcommand{\eps}{\varepsilon}
\newcommand{\ph}{\varphi}
\newcommand{\thet}{\vartheta}

\newcommand{\Aut}{\operatorname{Aut}}

\newcommand{\Conv}{\operatorname{Conv}}
\newcommand{\End}{\operatorname{End}}
\newcommand{\Koop}{\operatorname{Koop}}

\newcommand{\Iso}{\operatorname{Iso}}
\newcommand{\dist}{\operatorname{dist}}

\newcommand{\diam}{\operatorname{diam}}

\newcommand{\id}{\operatorname{id}}

\newcommand\qed{{\unskip\nobreak\hfil\penalty50
  \hskip2em\hbox{}\nobreak\hfil\mbox{\rule{1ex}{1ex} \qquad}
    \parfillskip=0pt \finalhyphendemerits=0\par\medskip}}

\begin{document}

\title
{$\mu$-norm and regularity}
\author{D.Treschev \\
Steklov Mathematical Institute of Russian Academy of Sciences
}
\date{}
\maketitle

\begin{abstract}
In \cite{Tre_PSI20} we introduce the concept of a $\mu$-norm for a bounded operator in a Hilbert space. The main motivation is the extension of the measure entropy to the case of quantum systems. In this paper we recall the basic results from \cite{Tre_PSI20} and present further results on the $\mu$-norm. More precisely, we specify three classes of unitary operators for which the $\mu$-norm generates a bistochastic operator. We plan to use the latter in the construction of quantum entropy.
\end{abstract}

\section{Introduction}
\label{sec:intro}

Let $\calX$ be a nonempty set and let $\calB$ be a $\sigma$-algebra of subsets $X\subset\calX$. Consider the measure space $(\calX,\calB,\mu)$, where $\mu$ is a probability measure: $\mu(\calX) = 1$.

Consider the Hilbert space $\calH = L^2(\calX,\mu)$ with the scalar product and the norm
$$
  \langle f,g\rangle = \int_\calX f\overline g \, d\mu, \quad
  \|f\| = \sqrt{\langle f,f\rangle} .
$$
For any bounded operator $W$ on $\calH$ let $\|W\|$ be its $L^2$ norm:
$$
  \|W\| = \sup_{\|f\|=1}  \|Wf\|.
$$

We say that $\chi = \{Y_1,\ldots,Y_J\}$ is a (finite, measurable) partition (of $\calX$) if
$$
  Y_j\in\calB, \quad
  \mu\big(\calX\setminus \cup_{1\le j\le J} Y_j\big) = 0, \quad
  \mu(Y_j\cap Y_k) = 0 \quad
  \mbox{ for any $j,k\in\{1,\ldots,J\}$, $k\ne j$}.
$$

We say that $\kappa = \{X_1,\ldots,X_K\}$ is a subpartition of $\chi = \{Y_1,\ldots,Y_J\}$ if for any $k\in\{1,\ldots,K\}$ there exists $j\in\{1,\ldots,J\}$ such that $\mu(X_k\setminus Y_j) = 0$.

For any $X\in\calB$ consider the orthogonal projector
\begin{equation}
\label{imath}
  \pi_X : \calH\to\calH, \qquad
  \calH\ni f \mapsto \pi_X f = \one_X \cdot f,
\end{equation}
where $\one_X$ is the indicator of $X$.

Let $W$ be a bounded operator on $\calH$. For any partition $\chi = \{Y_1,\ldots,Y_J\}$ we define
\begin{equation}
\label{calM}
  \calM_\chi(W) = \sum_{j=1}^J \mu(Y_j) \| W\pi_{Y_j} \|^2.
\end{equation}
In \cite{Tre_PSI20} we have introduced the definition of the $\mu$-norm:\footnote
{in fact, a seminorm}
\begin{equation}
\label{mumeasure}
  \|W\|_\mu = \inf_\chi \sqrt{\calM_\chi(W)}.
\end{equation}

Recall that the operator $U$ is said to be an isometry if
$$
  \langle f,g\rangle = \langle Uf,Ug\rangle, \qquad
  f,g\in\calH.
$$
If the isometry $U$ is invertible then $U$ is a unitary operator.

For any bounded $W$, any $Y\in\calB$, and any isometry $U$
$$
  \|W\pi_Y\| \le \|W\|, \quad
  \|UW\| = \|W\|, \quad
  \|\pi_Y\| = 1 \; \mbox{(if $\mu(Y)>0$)}.
$$
This implies the following obvious properties of the $\mu$-norm:
\begin{eqnarray}
\label{|1|}
         \|\id\|_\mu
     &=& 1,\qquad
         \|W\|_\mu
 \;\le\; \|W\|,     \\
\label{WW}
          \|W_1 W_2\|_\mu
    &\le& \|W_1\| \|W_2\|_\mu, \\
\label{lambdaW}
         \|\lambda W\|_\mu
     &=& |\lambda|\, \|W\|_\mu \quad
         \mbox{ for any $\lambda\in\mC$}, \\
\label{UWWU}
         \|W\|_\mu
     &=& \|UW\|_\mu \quad \mbox{ for any isometry $U$}.
\end{eqnarray}

The $\mu$-norm is motivated by the problem of the extension of the measure entropy\footnote
{known also as the Kolmogorov-Sinai entropy} to the case of quantum systems. Now we describe briefly the idea.

Let $F:\calX\to\calX$ be an endomorphism of the probability space $(\calX,\calB,\mu)$. This means that for any $X\in\calB$ the set $F^{-1}(X)$ (the complete preimage) also lies in $\calB$ and $\mu(X) = \mu(F^{-1}(X))$. Invertible endomorphisms are called automorphisms. Let $\End(\calX)$ denote the semigroup of all endomorphisms of $(\calX,\calB,\mu)$. There are two standard constructions.
\medskip

(1) Any $F\in\End(\calX,\mu)$ generates the isometry (a unitary operator if $F$ is an automotphism) $U_F$ on $\calH$ (the Koopman operator):
$$
  L^2(\calX,\mu)\ni f \mapsto U_F f = f\circ F, \qquad
  U_F =: \Koop(F).
$$

(2) For any $F\in\End(\calX,\mu)$ it is possible to compute the measure entropy $h(F)$.
\medskip

Our question is as follows. Is it possible to determine in some ``natural way'' a real nonnegative function $\frakh$ on the semigroup of isometries $\Iso(\calH)$ so that the diagram
$$
\begin{array}{rcl}
              &   \End(\calX,\mu)                  &      \\
  h\swarrow   &                                    &  \searrow\Koop \\
 \mR_+        & \stackrel{\frakh}{\longleftarrow}  &  \Iso(\calH)
\end{array}
$$
is commutative?

Recall the construction of the measure entropy of an endomorphism. Let $J_N$ be the set of multiindices $j = (j_0,\ldots,j_N)$, where any component $j_n$ takes values in the set $\{0,\ldots,K\}$. For any partition $\chi = \{X_0,\ldots,X_K\}$ and $j\in J_N$ we define
$$
  \bX_j = F^{-N}(X_{j_N})\cap \ldots \cap F^{-1}(X_{j_1}) \cap X_{j_0}.
$$
We define $h_F(\chi,N+1)$ by
$$
    h_F(\chi,N+1)
  = - \sum_{j\in J_N} \mu(\bX_j) \log \mu(\bX_j).
$$
The function $h_F$, as a function of the second argument, is subadditive:
$h_F(\chi,n+m) \le h_F(\chi,n) + h_F(\chi,m)$. This implies existence of the limit
$$
  h_F(\chi) = \lim_{n\to\infty} \frac1n h_F(\chi,n).
$$
Finally, the measure entropy is defined by
$$
  h(F) = \sup_\chi h_F(\chi).
$$

A rough idea is to construct the entropy of a unitary operator $U$ analogously with the following difference. Instead of $\bX_j$ we take
\begin{equation}
\label{frakX}
    \frakX_j
  = \pi_{X_{j_N}} U \pi_{X_{j_{N-1}}} U \ldots U \pi_{X_{j_1}} U \pi_{X_{j_0}} .
\end{equation}
We define
\begin{equation}
\label{hUchiN}
    \frakh_U(\chi,N+1)
  = - \sum_{j\in J_N} \|\frakX_j\|_\mu^2 \log \|\frakX_j\|_\mu^2 .
\end{equation}
Other details are the same:
\begin{equation}
\label{limfrakh}
    \frakh_U(\chi)
  = \lim_{n\to\infty} \frac1n \frakh(\chi,n), \qquad
    \frakh(U)
  = \sup_\chi \frakh_U(\chi)
\end{equation}
(provided the limit (\ref{limfrakh}) exists).

In \cite{Tre_PSI20} we prove that for any automorphism $F$
\begin{equation}
\label{UF=F}
  \frakh(U_F) = h(F) .
\end{equation}

In the literature there exist several attempts to extend the concept of the measure entropy to quantum systems, see \cite{CNT,Ohya95,Ohya00,AOW,AF,M,GLW} and many others. In \cite{A} several mutual relations between these approaches are given. Some works (for example, \cite{AOW, Sr, Pe, BG, KK}) deal with the finite-dimensional case $(\#\calX<\infty)$. In \cite{DF05, DF18} a construction for the measure entropy is proposed for doubly stochastic (bistochastic) operators on various spaces of functions on a measure space. It remains unclear, which approach to quantum generalization of the measure entropy is ``more physical''. This issue may become more clear after computation of the entropy $\frakh$ or its analogs in examples. We plan to do this in forthcoming papers.
\medskip

The above definition of $\frakh(U)$ meets several technical problems, including the question on subadditivity\footnote
{This subadditivity is important for the existence of the limit (\ref{limfrakh}).}
of $\frakh_U(\chi,n)$ and the inequality $\frakh_U(\chi)\le\frakh_U(\kappa)$ if $\kappa$ is a subpartition of $\chi$\footnote
{This inequality is necessary if we want to approach the supremum (\ref{limfrakh}) on fine partitions $\chi$.}.
We plan to change slightly the definition by replacing $\|\frakX_j\|_\mu^2$ in (\ref{hUchiN}) by another quantity similar to it to satisfy these two properties. To this end we associate with $U$ a bistochastic operator on $L^1(\calX,\mu)$. Then the entropy of $U$ may be defined analogously to \cite{DF05, DF18}. We will present details in another paper, but seems, this bistochastic operator may turn to be interesting by itself. As we will see below in this paper, the concept of the $\mu$-norm remains central in all our constructions.

We use the $\mu$-norm to construct a bistochastic operator in three cases: Koopman operators, the case $\#\calX < \infty$ and regular operators (a special class of operators defined in Section \ref{sec:reg_oper}) for $\calX=\mT$. It is more or less clear that the circle $\mT$ may be replaced by the torus $\mT^d$, but this will be a subject of another paper.

Before systematic attempts to compute entropy of various unitary operators and to study its proprties we have to study the $\mu$-norm $\|\cdot\|_\mu$. We have started this in \cite{Tre_PSI20}. The present paper is a continuation of this program.

\section{Notation and previous results}

Here we collect basic results from \cite{Tre_PSI20}. We will refer to some of them below.
\medskip

{\bf (1)}. $\|\pi_X\|_\mu^2 = \mu(X)$ for any $X\in\calB$.
\medskip

{\bf (2)}. If $\chi'$ is a subpartition of $\chi$ then $\calM_{\chi'}(W) \le \calM_\chi(W)$. Hence the quantities $\calM_\chi(W)$ approach the infinum (\ref{mumeasure}) on fine (small scale) partitions.
\medskip

{\bf (3)}. For any two bounded operators $W_1$ and $W_2$
$$
  \|W_1 + W_2\|_\mu \le \|W_1\|_\mu + \|W_2\|_\mu.
$$
This triangle inequality combined with (\ref{lambdaW}) imply that $\|\cdot\|_\mu$ is a seminorm on the space of bounded operators on $\calH$.
\medskip

{\bf (4)}. Let $F$ be an automorphism of $(\calX,\calB,\mu)$ and let $U_F=\Koop(F)$. Then
\begin{equation}
\label{UFpi}
  U_F\pi_X = \pi_{F^{-1}(X)} U_F \quad
  \mbox{for any } X\in\calB .
\end{equation}
For any bounded operator $W$
\begin{equation}
\label{WUF}
  \|W U_F\|_\mu = \|W\|_\mu.
\end{equation}
  This implies $\|U_F^{-1} W U_F\|_\mu = \|W\|_\mu$. Informally speaking, this means that measure preserving coordinate changes on $\calX$ preserve the $\mu$-norm.
\medskip

{\bf (5)}. $\|\cdot\|_\mu$ is a continuous function in the $L^2(\calX,\mu)$ operator topology.
\medskip

{\bf (6)}. If the measure $\mu$ has no atoms then $\|W+W_0\|_\mu=\|W\|_\mu$ for any bounded $W$ and compact $W_0$. In particular, $\mu$-norm of any compact operator vanishes. In fact, there exist non-compact operators with zero $\mu$-norm.
\medskip

{\bf (7)}. Given $g\in L^\infty(\calX,\mu)$ let $\widehat g$ be the multiplication operator defined by
$f\mapsto \widehat g f = gf$. Then $\|\widehat g\|_\mu = \|g\|$.
\medskip

{\bf (8)}. Suppose $\calX = \{1,\ldots,J\}$ is finite and the measure of any element equals $1/J$. Then $\calH$ is isomorphic to $\mC^J$ with the Hermitian product $\langle f,g\rangle_J = \frac1J\sum_{j=1}^J f(j)\overline{g(j)}$. Let
$$
  f\mapsto Wf, \quad  (Wf)(k) = \sum_{j=1}^J W(k,j) f(j)
$$
be an operator on $\calH$. Then
\begin{equation}
\label{|.|mu(finite)}
  \|W\|_\mu^2 = \frac1J \sum_{j,k=1}^J |W(k,j)|^2 .
\end{equation}
\medskip

{\bf (9)}. For any partition $\{X_1,\ldots,X_K\}$ of $\calX$
$$
     \|W\|_\mu^2
  =  \sum_{k=1}^K \|W \pi_{X_k}\|_\mu^2 , \quad
     \|W\|_\mu^2
 \le  \sum_{k=1}^K \|\pi_{X_k} W\|_\mu^2.
$$

{\bf (10)}. Let $\calX$ be a compact metric space and $\mu$ a Borel measure w.r.t. the corresponding topology. Let $B_r(x)\subset\calX$ denote the open ball with center at $x$ and radius $r$. Then for any $x\in\calX$ the limit
$$
  \thet(x) = \lim_{\eps\searrow 0} \|W\pi_{B_\eps(x)}\|^2
$$
exists, the function $\thet$ is measurable and $\|W\|_\mu^2 \le \int_\calX \thet \,d\mu$. There exists an example which shows that in general this inequality is strict. However $\|W\|_\mu^2 = \int_\calX \thet \,d\mu$ provided two additional conditions {\bf C1} and {\bf C2} hold:
\medskip

{\bf C1}. The function $\thet$ is continuous.

{\bf C2}. There exists $c>0$ such that for any open $\calO\subset\calX$, $\diam(\calO)\le\eps$ and any $x\in\calO$ there exists a function $f = \pi_\calO f$ satisfying
$$
      \Big| \|Wf\|^2 - \thet(x) \|f\|^2 \Big|
  \le \gamma(\eps) \|f\|^2, \quad
      c < \big| f|_\calO \big| < c^{-1},
$$
where $\gamma(\eps)\to 0$ as $\eps\to 0$. As usual, $f|_\calO$ denotes restriction of $f$ to the set $\calO$.
\medskip

The further results from \cite{Tre_PSI20} we mention here concern the case
$\calX = \mT = \mR / 2\pi\mZ$ with the Lebesgue measure $\mu$. We expect that they can be extended to the case $\calX = \mT^d$, $d>1$.
\medskip

{\bf (11)}. Let $\calX = \mT$ be a circle with the Lebesgue measure $d\mu = \frac1{2\pi} dx$. For any bounded sequence $\{\lambda_k\}_{k\in\mZ}$ we consider the distribution $\lambda(x) = \sum \lambda_k e^{ikx}$. The convolution operator
$$
      f\mapsto\Conv_\lambda f = \lambda * f := \int_\mT \lambda(y) f(\cdot - y)\, dy
$$
is bounded: $\|\Conv_\lambda\| = \sup_{k\in\mZ} |\lambda_k|$. Then
$$
      \|\Conv_\lambda\|_\mu^2 = \rho(\lambda), \qquad
      \rho(\lambda) = \limsup_{\# I\to\infty} \rho_I(\lambda), \quad
      \rho_I(\lambda) = \frac1{\# I} \sum_{k\in I} |\lambda_k|^2,
$$
where $I\subset\mZ$ are integer intervals.
\medskip

{\bf (12)}. One of the main technical tools in the analysis of the $\mu$-norm in the case $\calX = \mT$ is the following lemma on Fourier coefficients of localized functions on the circle.

\begin{lem}
\label{lem:YfJ}
Let $Y=[a-\eps,a+\eps]$ and
$\displaystyle   f = \pi_Y f = \sum_{k\in\mZ} f_k e^{ikx} \in L^2(\mT)$.
Then for any integer interval $J$ and any $m\in\mZ$
\begin{eqnarray}
\label{f-ef}
       \|f - e^{im(x-a)} f\|
 &\le& |m| \eps \|f\| , \\
\label{f-fe}
       |f_m - e^{ila} f_{m+l}|
 &\le& \frac{\eps^{3/2}}{\sqrt\pi} |l| \|f\| , \\
\label{ff-f}
       \Big| \sum_{k\in\mZ} e^{-ima} f_k\overline f_{k+m} - \|f\|^2 \Big|
 &\le& |m| \eps \|f\|^2 .
\end{eqnarray}
\end{lem}

{\bf (13)}. Consider the operator $W = (W_{j,k})_{j,k\in\mZ^d}$ on $\calH = L^2(\mT)$:
$$
  f = \sum_{k\in\mZ} f_k e^{ikx} \mapsto Wf = \sum_{j,k\in\mZ} W_{j,k} f_k e^{ijx}.
$$
We say that $W$ is of diagonal type ($W\in\DT(\mT)$) if
$$
  \sup_{j\in\mZ} |W_{j+k,j}| = c_k < \infty, \quad k\in\mZ\quad
  \mbox{and}\quad \sum_{k\in\mZ} c_k = \bc < \infty.
$$
The sequence $c_s$ is said to be the majorating sequence for $W\in\DT(\mT)$. We define the norm $\|W\|_\DT = \bc$.

Operators from $\DT(\mT)$ are bounded. As simple examples we have the following operators of diagonal type.
\smallskip

(a) Bounded convolution operators.

(b) Operators of multiplication by functions with absolutely converging Fourier series.

(c) The conjugated operator $W^*$ if $W\in\DT(\mT)$.

(d) Linear combinations and products of operators of diagonal type. Moreover,
$$
  \|W' W''\|_\DT \le \|W'\|_\DT \|W''\|_\DT \quad
  \mbox{for all $W',W''\in\DT(\mT)$}.
$$

{\bf (14)}. We prove that the normed space $\big(\DT(\mT), \|\cdot\|_\DT\big)$ is closed. As a corollary we obtain that $\big(\DT(\mT), \|\cdot\|_\DT\big)$ is a $C^*$-algebra.
\medskip

{\bf (15)}. We have the following inequalities between the norms:
\begin{eqnarray}
\label{.<.<.}
&  \|W\|_\mu \le \|W\| \le \|W\|_\DT  \quad
   \mbox{for any $W\in\DT(\mT)$} , & \\
\label{||f||<||f||_DT}
&  \|f\|_\infty \le \|\widehat f\|_\DT \quad
    \mbox{for any $f$ with absolutely converging Fourier series} ,  &
\end{eqnarray}
where $\widehat f$ is the operator of multiplication by $f$.

{\bf (16)}. We associate with $W\in\DT(\mT)$ and any point $a\in\mT$ the distribution $L_a$,
\begin{equation}
\label{L_a}
  L_a = \sum_{j\in\mZ} w_j(a) e^{ijx}, \qquad
  w_j(a) = \sum_{k\in\mZ} W_{j,k} e^{i(j-k)a}.
\end{equation}
For any $l\in\mZ$ and $a\in\mT$ we have the estimate
\begin{equation}
\label{|w(a)|<}
   |w_l(a)| \le \bc = \|W\|_\DT .
\end{equation}

We prove that for any $W\in\DT(\mT)$ the function
$a\mapsto\rho(L_a) = \limsup_{\# I\to\infty} \rho_I(L_a)$ ($I\subset\mZ$ are intervals, $\rho_I$ is defined in item {\bf (11)}) is continuous and
\begin{equation}
\label{||W||_mu=int}
     \|W\|_\mu^2 = \frac1{2\pi} \int_\mT \rho(L_a) \, da .
\end{equation}

{\bf (17)}. For any operator $W\in\DT(\mT)$ we introduce the average trace of $W^* W$ by
$$
    \bT(W)
  = \limsup_{\# I\to\infty} \frac1{\# I} \sum_{j\in\mZ,\, l\in I} |W_{l,j}|^2,
$$
where $I\subset\mZ$ are intervals. Then
\begin{equation}
\label{||W||_mu>}
    \bT(W) \le \|W\|_\mu^2.
\end{equation}

{\bf (18)}. We also prove that if $U\in\DT(\mT)$ is a unitary operator then
\begin{equation}
\label{T=T=T}
  \bT(W) = \bT(WU) = \bT(UW).
\end{equation}

\section{Main results}

In this section we collect main results of the present paper. In short, these results concern (1) properties of the $\mu$-norm on $\calR(\mT)$, a special class of operators on $L^2(\mT)$, and (2) a construction of a bistochastic operator $\calW$ on $L^1(\calX,\mu)$ associated with an operator $W$ under some conditions, imposed on $W$.

\begin{itemize}

\item We say that an operator $W\in\DT(\mT)$ is regular (the notation is $W\in\calR(\mT)$) if for any $m,n\in\mZ$ there exists the limit
$$
     \omega_{m,n}
  =  \lim_{\# I\to\infty} \omega_{I,m,n} ,  \qquad
     \omega_{I,m,n}
  =  \frac1{\# I} \sum_{j\in\mZ,\, l\in I}
        W_{l+m,j} \overline W_{l,j+n} .
$$

\item By Lemma \ref{lem:Rclosed} $\calR(\mT)$ is a cone closed with respect to the norm $\|\cdot\|_\DT$.

\item By Lemma \ref{lem:dimregular} $\|W\|_\mu^2 = \bT(W)$ for any $W\in\calR(\mT)$ (compare with (\ref{||W||_mu>})).

\item Let $\ACF(\mT)$ be the space of functions on $\mT$ with absolutely converging Fourier series. Then the space of multiplication operators $\widehat g$ by functions $g\in\ACF(\mT)$ form a $C^*$-subalgebra in the $C^*$-algebra $\DT(\mT)$.

\item In Section \ref{sec:reg_exa} we present several examples of regular operators.

    Let $W\in\DT(\mT)$ be an operator with periodic matrix i.e., there exists $\tau\in\mN$ such that $W_{j+\tau,k+\tau} = W_{j,k}$ for any $j,k\in\mZ$. Any such operator is regular (Lemma \ref{lem:per=>reg}).

  Another example of a regular operator is $\Conv_\lambda$, where
  $\lambda_k=e^{i\tau k^2}$.

\item By Proposition \ref{prop:WgWWg_reg} for any $W\in\calR(\mT)$ and $g_1,g_2\in\ACF(\mT)$ the operator $\widehat g_1 W \widehat g_2$ is regular.

\item In Section \ref{sec:mu_W} we associate with operators $W$ from the following three classes

  (1) Koopman operators $\Koop(F)$, $F\in\Aut(\calX,\mu)$,

  (2) operators in the case $\#\calX < \infty$,

  (3) operators from $\calR(\mT)$

  a measure $d\mu_W(x',x'') = \nu(x',x'')\, d\mu(x') d\mu(x'')$ on $\calX\times\calX$ such that for any ``sufficiently regular'' functions $g',g'' : \calX\to\mC$
  $$
    \|\widehat g' W\widehat g''\|_\mu^2 = \int_{\calX^2} |g'(x')|^2 |g''(x'')|^2 \, d\mu_W(x',x'') .
  $$

\item In Section \ref{sec:bistochastic} we introduce the operator
$$
  L^1(\calX,\mu) \ni f \mapsto \calW f = \int_\calX \nu(\cdot, a) f(a)\, d\mu(a)
$$
and prove (Lemma \ref{lem:bistochastic}) that $\calW$ is a bistochastic operator on $L^1(\calX,\mu)$.

\end{itemize}

\section{Regular operators}
\label{sec:reg_oper}

\subsection{Definition of $\omega_{m,n}$}

\begin{dfn}
\label{dfn:regular}
We say that $W\in\DT(\mT)$ is regular ($W\in\calR(\mT)$) if for any $m,n\in\mZ$ there exists the limit
\begin{equation}
\label{limomega}
    \lim_{\# I\to\infty} \omega_{I,m,n}
  = \omega_{m,n}, \qquad
    \omega_{I,m,n}
  = \frac1{\# I} \sum_{j\in\mZ,\, l\in I} W_{l+m,j} \overline W_{l,j+n} ,
\end{equation}
where $I$ are integer intervals.
\end{dfn}

Note that $\omega_{0,0}$ coincides with the average trace of $W^* W$:
\begin{equation}
\label{omega=T}
  \omega_{0,0} = \bT(W)\quad
  \mbox{if $W\in\calR(\mT)$}.
\end{equation}

For any integer interval $I$ we put
\begin{equation}
\label{vIm}
    v_{I,m}(a)
  = \sum_{l\in I} \frac{w_{l+m}(a) \overline w_l(a)}{\# I} ,
\end{equation}
where the functions $w_l(a)$ are defined in (\ref{L_a}).

\begin{lem}
\label{lem:lim_v}
Suppose $W$ is regular. Then for any $m\in\mZ$ there exists the limit
\begin{equation}
\label{limv}
    \lim_{\# I\to\infty} v_{I,m}(a)
  = v_m(a), \quad
    v_m(a) = \sum_{n\in\mZ} \omega_{m,n} e^{i(m+n)a}
\end{equation}
uniformly in $a\in\mT$. The Fourier series of the function $v_m$ absolutely converges.
\end{lem}

{\it Proof}. By (\ref{L_a})
\begin{eqnarray*}
     v_m(a)
 &=& \lim_{\# I\to\infty} \frac1{\# I} \sum_{j,k\in\mZ,\,l\in I}
           W_{l+m,j} \overline W_{l,k} e^{i(m-j+k)a} \\
 &=& \lim_{\# I\to\infty} \sum_{n\in\mZ}
           \omega_{I,m,n} e^{i(m+n)a}
 \; =\;  \sum_{n\in\mZ} \omega_{m,n} e^{i(m+n)a} .
\end{eqnarray*}
By (\ref{omegaI}) the limit is uniform in $a$.
By (\ref{sumomega}) this Fourier series absolutely converges. \qed

Note that by (\ref{|w(a)|<}) for any interval $I\subset\mZ$, any $m\in\mZ$, and any $a\in\mT$
\begin{equation}
\label{|v|}
  |v_{I,m}(a)| \le \bc^2, \quad
     |v_m(a)|  \le \bc^2.
\end{equation}

\subsection{Closeness with respect to $\|\cdot\|_{\DT}$}

If $W\in\calR(\mT)$ then for any $\lambda\in\mC$ the operator $\lambda W$ is also regular. Hence, regular operators form a cone $\calR(\mT)\subset\DT(\mT)$.

\begin{lem}
\label{lem:Rclosed}
The cone $\calR(\mT)$ is closed with respect to the norm $\|\cdot\|_{\DT}$.
\end{lem}

{\it Proof}. Suppose $\{W_p\}_{p\in\mN}$, $W_p\in\calR(\mT)$ is a Cauchy sequence. By {\bf (14)} there exists $W = \lim_{p\to\infty} W_p$, where the limit is taken with respect to the norm $\|\cdot\|_{\DT}$.

For any $\eps>0$ there exists positive $N$ such that
\begin{equation}
\label{Wp-Wq}
   \mbox{for any integer $p,q>N$ we have:} \quad
   \|W_p - W_q\|_\DT < \eps.
\end{equation}
We define $(\omega_p)_{I,m,n}$ and $(\omega_p)_{m,n}$ by (\ref{limomega}), where $W$ is replaced by $W_p$. Then
\begin{eqnarray*}
     \big| (\omega_p)_{I,m,n} - (\omega_q)_{I,m,n} \big|
 &=& \frac1{\# I}
      \bigg| \sum_{j\in\mZ, l\in I} \bigg(
             (W_p)_{l+m,j} (\overline W_p)_{l,j+n}
           - (W_q)_{l+m,j} (\overline W_q)_{l,j+n}
      \bigg)\bigg|   \\
 &\le& \frac1{\# I} \big( \Sigma_1 + \Sigma_2 \big), \\
      \Sigma_1
 &=&  \sum_{j\in\mZ, l\in I} \bigg|
         (W_p)_{l+m,j} \bigg( (\overline W_p)_{l,j+n} - (\overline W_q)_{l,j+n}
                       \bigg) \bigg|, \\
      \Sigma_2
 &=&  \sum_{j\in\mZ, l\in I} \bigg|
         \bigg( (W_p)_{l+m,j} - (W_q)_{l+m,j} \bigg) (\overline W_q)_{l,j+n}
                        \bigg| .
\end{eqnarray*}

To estimate the sums $\Sigma_1$ and $\Sigma_2$, we put
$$
  (c_p)_k = \sup_{j\in\mZ} |(W_p)_{k+j,j}|, \quad
  (c_q)_k = \sup_{j\in\mZ} |(W_q)_{k+j,j}|, \quad
   d_k    = \sup_{j\in\mZ} |(W_p)_{k+j,j} - (W_q)_{k+j,j}|.
$$
The sums $\sum_k (c_p)_k$ are uniformly bounded:
\begin{equation}
\label{cpcq}
  \sum_{k\in\mZ} (c_p)_k \le \tilde\bc, \quad
  \sum_{k\in\mZ} (c_q)_k \le \tilde\bc \quad
  \mbox{for some constant $\tilde\bc$}.
\end{equation}
Moreover, by (\ref{Wp-Wq})
\begin{equation}
\label{dk}
  \sum_{k\in\mZ} d_k < \eps.
\end{equation}

By (\ref{cpcq}) and (\ref{dk})
$$
      \Sigma_1
  \le \sum_{j\in\mZ,\, l\in I} (c_p)_{l+m-j} d_{l-j-n}
  \le \# I \tilde\bc \eps.
$$
Analogously $\Sigma_2 \le \# I \tilde\bc \eps$. This implies that for any interval $I$
$$
  \big| (\omega_p)_{I,m,n} - (\omega_q)_{I,m,n} \big| \le 2\tilde\bc \eps.
$$
Hence, $\big| (\omega_p)_{m,n} - (\omega_q)_{m,n} \big| \le 2\tilde\bc\eps$ i.e., for any integer $m,n$ the sequence $(\omega_p)_{m,n}$, $p\in\mN$ is a Cauchy sequence. \qed

\subsection{$\mu$-norm of a regular operator}

\begin{lem}
\label{lem:dimregular}
If $W\in\calR(\mT)$ then (compare with (\ref{||W||_mu>}))
\begin{equation}
\label{dim=omega}
    \|W\|_\mu^2
  = \bT(W) .
\end{equation}
\end{lem}

{\it Proof}. By (\ref{L_a}), (\ref{||W||_mu=int}) and (\ref{vIm})
\begin{equation}
\label{dim=intlimreg}
     \|W\|_\mu^2
  =  \frac1{2\pi} \int_\mT \lim_{\# I\to\infty} v_{I,0}(a) \, da.
\end{equation}
By Lemma \ref{lem:lim_v} the limit $v_0(a) = \lim_{\#I\to\infty} v_{I,0}(a)$ exists for any $a\in\mT$ and by (\ref{|v|}) $|v_{I,0}(a)| \le \bc^2$ for all $I$ and $a$. Therefore by the Lebesgue theorem on bounded convergence we may exchange the integration and the limit:
\begin{eqnarray*}
     \|W\|_\mu^2
 &=& \lim_{\# I\to\infty} \frac1{2\pi} \int_\mT v_{I,0}(a)\, da
  =  \lim_{\# I\to\infty} \frac1{2\pi} \int_\mT \sum_{j,k\in\mZ,\, l\in I}
         \frac1{\# I} W_{l,j} \overline W_{l,k} e^{i(k-j)a} \, da \\
 &=& \lim_{\# I\to\infty} \sum_{j\in\mZ,\, l\in I}
         \frac1{\# I} W_{l,j} \overline W_{l,j}
 \; =\; \omega_{0,0} .
\end{eqnarray*}
By (\ref{limomega}) this implies (\ref{dim=omega}). \qed

Assertion {\bf (18)} and Lemma \ref{lem:dimregular} imply the following.

\begin{cor}
\label{cor:UWU}
Suppose $W,U\in\DT(\mT)$, where $U$ is unitary and both $W$ and $UWU^{-1}$ are regular. Then by {\bf (18)} and (\ref{UWWU})
$$
  \|W\|_\mu = \|UWU^{-1}\|_\mu = \|WU^{-1}\|_\mu .
$$
\end{cor}

\section{Regular operators: examples}
\label{sec:reg_exa}

\begin{dfn}
We say that the matrix $(W_{k,j})$ of the operator $W\in\DT(\mT)$ is $\tau$-periodic, $\tau\in\mN$, if
$$
  W_{k+\tau,j+\tau} = W_{k,j} \quad\mbox{for any } k,j\in\mZ.
$$
\end{dfn}

In particular, for any $g\in\ACF(\mT)$ matrix of the operator $\widehat g\in\DT(\mT)$ is 1-periodic.

\begin{lem}
\label{lem:per=>reg}
Suppose $W\in\DT(\mT)$ is an operator with $\tau$-periodic matrix, $\tau\in\mN$. Then $W\in\calR(\mT)$ and $\omega_{m,n} = \breve\omega_{m,n}$,
$$
  \breve\omega_{m,n} = \frac1\tau \sum_{j\in\mZ,\,l\in J}
         W_{l+m,j} \overline W_{l,j+n},
$$
where $J\subset\mZ$ is any interval with $\# J = \tau$.
\end{lem}

{\it Proof}. By $\tau$-periodicity $\breve\omega_{m,n}$ does not depend on $J$.
For any $I = \{s,s+1,\ldots,s+K\}$ let $I'\subset I$ be the maximal integer interval of the form $I' = \{s,s+1,\ldots,s+q\tau - 1\}$, $q\in\mZ$. We put
$I'' = \{s+q\tau,s+q\tau+1,\ldots,s+K\}$. Then $\# I'' < \tau$,
$\omega_{I',m,n} = \breve\omega_{m,n}$ while
$\omega_{I,m,n} - \omega_{I',m,n} = A_1 + A_2$,
$$
     A_1
  =  \frac1{\# I} \sum_{l\in I''\setminus I'}
                   W_{l+m,j} \overline W_{l,j+n}, \quad
     A_2
  =  \Big(\frac1{\# I} - \frac1{\# I'}\Big) \sum_{l\in I'}
                   W_{l+m,j} \overline W_{l,j+n}.
$$
The inequalities
$$
    |A_1| \le \frac{\tau\bc^2}{\# I}, \quad
    |A_2| \le \frac{\tau\bc^2}{\# I - \tau}
$$
imply the existence of the limit (\ref{limomega}) and the equation
$\omega_{m,n} = \breve\omega_{m,n}$. \qed

\begin{cor}
Consider the operator $\widehat g$, where $g\in\ACF(\mT)$. After simple calculations we obtain:
\begin{equation}
\label{omega(g)}
  \omega_{m,n}(\widehat g) = \bg_{m+n}, \qquad
  g(x) \overline g(x) = \sum_{k\in\mZ} \bg_k e^{ikx},
\end{equation}
where the last equation is the definition of $\bg_k$.
\end{cor}

\begin{lem}
\label{lem:su}
Let $\{\lambda_k\}_{k\in\mZ}$ be defined by $\lambda_k = e^{i\tau k^2}$. Then the operator $\Conv_\lambda$ is regular:
$$
    \omega_{m,n}
  = \delta_{0,m+n} \delta_{\tau m,\pi\mZ} e^{i\tau m^2}, \qquad
    \delta_{\tau m,\pi\mZ}
  = \left\{ \begin{array}{cl}
             1 & \mbox{if } \tau m/\pi \in \mZ , \\
             0 & \mbox{if } \tau m/\pi \not\in \mZ .
            \end{array}
    \right.
$$
\end{lem}

{\it Proof}. In this case $W_{k,j} = \delta_{kj} e^{i\tau k^2}$ and
\begin{eqnarray*}
     \omega_{I,m,n}
 &=& \frac1{\# I} \sum_{j\in\mZ,\,l\in I}
             \delta_{l+m,j}\delta_{l,j+n} e^{i\tau((l+m)^2 - l^2} \\
 &=& \frac1{\# I} \sum_{l\in I}
             \delta_{l,l+m+n} e^{i\tau m^2 + 2i\tau ml}
 \; =\; \frac1{\# I} \delta_{0,m+n} e^{i\tau m^2}
                     \sum_{l\in I} e^{2i\tau ml} .
\end{eqnarray*}
Hence, if $\tau m / \pi \in \mZ$ then $\omega_{m,n} = \delta_{0,m+n} e^{i\tau m^2}$. If $\tau m / \pi \not\in \mZ$ then $\omega_{m,n} = 0$. \qed

\section{Measure associated with an operator}
\label{sec:mu_W}

\subsection{Koopman operator}

Let $(\calX,\calB,\mu)$ be a probability space.

\begin{lem}
\label{lem:gU_F}
$F\in\Aut(\calX,\mu)$, $g_0,\ldots,g_K\in L^\infty(\calX,\mu)$. Then
\begin{equation}
\label{gU_F}
    \big\| \widehat g_K U_F \widehat  g_{K-1} \ldots U_F \widehat  g_0 \big\|^2_\mu
  = \int_\calX  |g_K\circ F^K|^2 |g_{K-1}\circ F^{K-1}|^2 \ldots |g_0|^2 \, d\mu
\end{equation}
\end{lem}

{\it Proof}. Take a small constant $\sigma>0$ and consider the partition $\{X_1,\ldots,X_J\}$ such that\footnote
{We can use the same partition for all the functions $g_k$.}
\begin{equation}
\label{g-Sigma}
  \| g_k - \ph_k \|_\infty < \sigma \|g_k\|_\infty ,
  \quad   \ph_k = \sum_{j=1}^J g_{k,j} \one_{X_j} ,
  \quad   \|\ph_k\|_\infty \le \|g_k\|_\infty ,
  \qquad  k = 0,\ldots,K .
\end{equation}
Here $g_{k,j}\in\mC$ are some constants.
We put $S = \widehat\ph_K U_F \widehat\ph_{K-1} \ldots U_F \widehat \ph_0$.
By the triangle inequality the quantity
$$
  \Big| \| \widehat g_K U_F \widehat g_{K-1} \ldots U_F \widehat g_0 \|_\mu
            - \|S\|_\mu \Big|
$$
does not exceed
\begin{eqnarray*}
  &&  \Big| \| \widehat g_K U_F \widehat g_{K-1} \ldots U_F \widehat g_0 \|_\mu
        - \| \widehat\ph_K U_F \widehat g_{K-1} \ldots U_F \widehat g_0 \|_\mu \Big|\\
 &+&  \Big| \| \widehat\ph_K U_F \widehat g_{K-1} \ldots U_F \widehat g_0 \|_\mu
        - \| \widehat\ph_K U_F \widehat g_{K-1} \ldots U_F \widehat g_0 \|_\mu \Big|\\
 &+&   \ldots
\; +\;  \Big| \| \widehat\ph_K U_F \widehat\ph_{K-1}
                         \ldots \widehat\ph_1 U_F \widehat g_0 \|_\mu
        - \| \widehat\ph_K U_F \widehat\ph_{K-1}
                         \ldots \widehat\ph_1 U_F \widehat\ph_0 \|_\mu \Big| \\
 \le && \| (\widehat g_K - \widehat\ph_K) U_F \widehat g_{K-1}
                      \ldots U_F \widehat g_0 \|_\mu
  +      \| \widehat\ph_K U_F (\widehat g_{K-1} - \widehat\ph_{K-1})
                        \ldots U_F \widehat g_0 \|_\mu \\
 &+&   \ldots \;
  +    \| \widehat\ph_K U_F \widehat\ph_{K-1}
                 \ldots \widehat\ph_1 U_F (\widehat g_0 - \widehat\ph_0) \|_\mu \\
 \le && (K+1) \sigma \prod_{0\le k\le K} \|g_k\|_\infty .
\end{eqnarray*}
By (\ref{UFpi})
$$
    S
 =  \sum_{j_0,\ldots,j_K} g_{K,j_K} g_{K-1,j_{K-1}} \ldots g_{0,j_0}
           \pi_{X_{j_K}} \pi_{F^{-1}(X_{j_{K-1}})} \ldots \pi_{F^{-K}(X_{j_0})} U_F^K.
$$
We put $X_{j_K,\ldots,j_0} = X_{j_K} \cap F^{-1}(X_{j_{K-1}}) \cap
                                           \ldots \cap F^{-K}(X_{j_0})$.
Note that for any two sets $X'=X_{j'_K,\ldots,j'_0}$ and $X''=X_{j''_K,\ldots,j''_0}$ we have: $\mu(X'\cap X'')=0$ if the collections of indices $j'_K,\ldots,j'_0$ and $j''_K,\ldots,j''_0$ do not coincide.
Then by (\ref{WUF}) and {\bf (1)} we obtain:
\begin{eqnarray}
\nonumber
           S
 &=& \sum_{j_0,\ldots,j_K} | g_{K,j_K} g_{K-1,j_{K-1}} \ldots g_{0,j_0} |^2
               \| \pi_{X_{j_K,\ldots,j_0}} \|_\mu^2 \\
\label{gU_FgU_F}
 &=&     \sum_{j_0,\ldots,j_K} | g_{K,j_K} g_{K-1,j_{K-1}} \ldots g_{0,j_0} |^2
            \mu( X_{j_K,\ldots,j_0} ) .
\end{eqnarray}

Equation (\ref{g-Sigma}) implies
$$
      \ph_k\circ F^k
  =   \sum_{j=1}^J g_{k,j} \one_{F^{-k}(X_j)} .
$$
Hence (\ref{gU_FgU_F}) is an integral sum for the integral
$$
  \int_\calX  |\ph_K\circ F^K|^2 |\ph_{K-1}\circ F^{K-1}|^2 \ldots |\ph_0|^2 \, d\mu .
$$
This integral differs from the integral (\ref{gU_F}) at most by
$2\sigma\prod_{0\le k\le K} \|g_k\|_\infty^2$. Since $\sigma$ is arbitrarily small, we obtain equation (\ref{gU_F}). \qed
\smallskip

Now suppose that $\calX$ is in addition a topological space and $\calB$ is the corresponding Borel $\sigma$-algebra.
Consider the distribution $\delta \in (C(\calX^2))^*$ (a measure on $\calX^2$ by the Riesz theorem) such that for any $\Phi\in C(\calX^2)$
\begin{equation}
\label{delta}
    \int_{\calX^2} \delta(x',x'') \Phi(x',x'')\, d\mu(x') d\mu(x'')
  = \int_\calX \Phi(x,x)\, d\mu(x) .
\end{equation}
Taking in (\ref{delta}) $\Phi(x',x'')=\ph'(x')\ph''(x'')$, where
$\ph',\ph'' : \calX\to\mC$ are arbitrary continuous functions, we obtain:
\begin{equation}
\label{int_delta}
  \int_\calX \delta(x',x'') \ph'(x')\, d\mu(x') = \ph'(x''), \quad
  \int_\calX \delta(x',x'') \ph''(x'')\, d\mu(x'') = \ph''(x').
\end{equation}

Recall also that if $F\in\End(\calX,\mu)$
\begin{equation}
\label{intfF}
  \int_\calX f\, d\mu = \int_\calX f\circ F\, d\mu \quad
  \mbox{for any }  f\in L^1(\calX,\mu).
\end{equation}
Then for any continuous $F\in\Aut(\calX,\mu)$
\begin{equation}
\label{delta=delta}
   \delta(x',x'') = \delta(F(x'),F(x'')).
\end{equation}
Indeed, for any $\Phi\in C(\calX^2)$ by (\ref{intfF})
\begin{eqnarray*}
 &&    \int_{\calX^2} \delta(F(x'),F(x'')) \Phi(x',x'')\, d\mu(x') d\mu(x'') \\
 &=& \int_{\calX^2} \delta(x',x'') \Phi(F^{-1}(x'),F^{-1}(x''))\, d\mu(x') d\mu(x'')
 \; =\; \int_\calX \Phi(x,x)\, d\mu(x) .
\end{eqnarray*}

\begin{lem}
\label{lem:piUFpi}
Let $F\in\Aut(\calX,\mu)$ be continuous and $g_0,\ldots,g_K\in C(\calX)$. Then
\begin{eqnarray}
\!\!\!\!\!\nonumber
 &&    \big\|\widehat g_K U_F \ldots\widehat g_1 U_F \widehat g_0 \big\|_\mu^2 \\
\!\!\!\!\!\label{gU_F=int}
 &=&  \int_{\calX^{K+1}}  |g_K(x_K)|^2 \delta(x_K,F(x_{K-1})) \ldots
                          |g_1(x_1)|^2 \delta(x_1,F(x_0)) |g_0(x_0)|^2 \,
                             d\mu^{K+1} , \\
\!\!\!\!\!\nonumber
 && \qquad   d\mu^{K+1}
         =   d\mu(x_K)\ldots d\mu(x_0)  .
\end{eqnarray}
\end{lem}

{\it Proof}. Let $I$ denote the integral (\ref{gU_F=int}). We use in this integral the change of coordinates
$$
  x_0 = x'_0, \quad x_1 = F(x'_1), \quad \ldots \quad x_K = F^K(x'_K), \quad
  d\mu^{\prime K+1} = d\mu(x'_K)\ldots d\mu(x'_0)  .
$$
Then
\begin{eqnarray*}
     I
 &=& \int |g_K\circ F^K(x'_K)|^2 \delta(F^K(x'_K),F^K(x'_{K-1})) \\
 && \qquad\qquad
                  \ldots \; |g_1\circ F(x'_1)|^2 \delta(F(x'_1),F(x'_0)) |g_0(x'_0)|^2
                    \, d\mu^{\prime K+1}  \\
 &=&  \int_\calX |g_K\circ F^K|^2 \ldots |g_1\circ F|^2 |g_0|^2\, d\mu .
\end{eqnarray*}
It remains to use Lemma \ref{lem:gU_F}.  \qed
\smallskip

We associate with any $U_F=\Koop(F)$, where $F\in\Aut(\calX,\mu)$ is continuous,\footnote{We expect that the continuity is inessential here.} the measure $\mu_{U_F}$ on $\calX^2$:
\begin{equation}
\label{muUF}
  d\mu_{U_F}(x',x'') = \delta(x',F(x''))\, d\mu(x') d\mu(x'') .
\end{equation}

\subsection{The finite-dimensional case}
\label{sec:mu_W_finite}

Let $\calX = \{1,\ldots,J\}$ be a finite set. We identify $\calX$ with $\mZ_J = \mZ / J\mZ$, the cyclic additive group with $J$ elements.\footnote{Below we also use the structure of a commutative ring on $\mZ_J$.}
Let the measure $\mu$ of any point be equal to $1/J$. Then
$\calH = L^2(\calX,\mu)\cong (\mC^J,\langle\,,\rangle)$, where $\langle\,,\rangle$ equals the standard Hermitian product divided by $J$.

We put $\eta=e^{2\pi i/J}$. Then $\eta^J=1$ and $\overline\eta=\eta^{-1}$. The space $\calH$ may be identified with the space of ``discrete trigonometric polynomials''
\begin{equation}
\label{f(eta)}
  f = f(x) = \sum_{j\in\mZ_J} f_j \eta^{jx}, \qquad
  f_j\in\mC, \quad x\in\mZ_J .
\end{equation}
This polynomial representation generates on $\calH$ an operation of multiplication: for any two vectors
$f' = \sum f'_j \eta^{jx}$ and $f'' = \sum f''_j \eta^{jx}$
$$
  f' f'' = f = \sum_{k\in\mZ_J} f_k \eta^{kx}, \qquad
  f_k = \sum_{j\in\mZ_J} f'_{k-j} f''_j .
$$
This product introduces on $\calH$ the structure of a commutative ring. The structure of a Hilbert space is determined by
$$
  \langle f,g\rangle = \frac1J \sum_{k\in\mZ_J} f(k) \overline{g(k)} .
$$

The coefficients $f_k$ and $f(k)$ are connected by the ``discrete Fourier transform''
\begin{equation}
\label{discreteFourier}
  f(k) = \sum_{j\in\mZ_J} f_j \eta^{kj}, \quad
  f_j = \frac1J \sum_{k\in\mZ_J} f(k) \eta^{-kj} .
\end{equation}

For any $f$ satisfying (\ref{f(eta)}) and any $A\subset\mZ_J$ we put
$$
     \int_A f(x)\, d\mu(x)
  := \frac1J \sum_{k\in A} f(k)
   = \frac1J \sum_{k\in A,\, j\in\mZ_J} f_j \eta^{jk} .
$$
Then
$$
  \int_\calX f\, d\mu = f_0, \quad
  \langle f,g\rangle = \frac1J \int_\calX f\overline g \, d\mu .
$$

Consider a linear operator
$$
  f\mapsto Wf, \qquad
  (Wf)_k = \sum_{j\in\mZ_J} W_{kj} f_j .
$$
In another basis it takes the form
$$
  (Wf)(k) = \sum_{j\in\mZ_J} W(k,j) f(j) .
$$
Equations (\ref{discreteFourier}) imply
\begin{equation}
\label{WfWf}
  W(m,n) = \frac1J \sum_{j,k\in\mZ_J} \eta^{mk} W_{kj} \eta^{-jn}.
\end{equation}

We define
\begin{equation}
\label{nu_finite}
     \omega_{m,n}
  =  \frac1J \sum_{j,l\in\mZ_J} W_{l+m,j} \overline W_{l,j+n} , \quad
     \nu(x,a)
  =  \sum_{m,n\in\mZ_J} \omega_{m,n} \eta^{mx} \eta^{na} , \qquad
     x,a\in\mZ_J.
\end{equation}

\begin{lem}
\label{lem:WW=WW}
$\sum_{m,n\in\mZ_J} |W(m,n)|^2 = \sum_{j,k\in\mZ_J} |W_{j,k}|^2$.
\end{lem}

{\it Proof}. Direct computation with the help of (\ref{WfWf}) and the identity $\sum_{j\in\mZ_J}\eta^{jk}=J\delta_{k,0}$.  \qed

\begin{cor}
By (\ref{|.|mu(finite)}) and Lemma \ref{lem:WW=WW}
\begin{equation}
\label{||mu=omega00}
    \|W\|_\mu^2
  = \omega_{0,0}
  = \int_{\calX^2} \nu(x,a)\, d\mu(x) d\mu(a) .
\end{equation}
\end{cor}

\begin{lem}
\label{lem:gnug}
For any $g',g''\in\calH$ the operator $\widetilde W = \widehat g' W\widehat g''$ generates the coefficients $\widetilde\omega_{m,n}$ such that
\begin{equation}
\label{g'nug''}
    \sum_{m,n\in\mZ} \widetilde\omega_{m,n} \eta^{mx} \eta^{na}
  = |g'(x)|^2 \nu(x,a) |g''(a)|^2 .
\end{equation}
The measure
$d\mu_W(x,a) = \nu(x,a)\, d\mu(x) d\mu(a)$ on $\calX^2$ satisfies
\begin{equation}
\label{g'g''}
    \|\widehat g' W \widehat g''\|_\mu^2
  = \int_{\calX^2} |g'(x)|^2 |g''(a)|^2 \, d\mu_W(x,a) .
\end{equation}
\end{lem}

{\it Proof}. We put
$$
    g'(x) = \sum g'_k\eta^{kx}, \quad
    g''(x) = \sum g''_k\eta^{kx}, \quad
  |g'(x)|^2 = \sum \bg'_k\eta^{kx}, \quad
  |g''(x)|^2 = \sum \bg''_k\eta^{kx} .
$$
The equation $\widetilde W_{p,q} = \sum_{\alpha,\beta} g'_{p-\alpha} W_{\alpha,\beta} g''_{\beta-q}$ implies
\begin{eqnarray*}
     \widetilde\omega_{m,n}
 &=& \frac1J \sum_{l,j,\alpha,\beta} g'_{l+m-\alpha'} W_{\alpha',\beta'} g''_{\beta'-j}
                   \overline g'_{l-\alpha''} \overline W_{\alpha'',\beta''} \overline g''_{\beta''-j-n} \\
 &=& \frac1J \sum_{\alpha,\beta} \bg'_{m-\alpha'+\alpha''} W_{\alpha',\beta'}
                                  \overline W_{\alpha'',\beta''} \bg''_{\beta'-\beta''+n} \\
 &=& \frac1J \sum_{p,q,\alpha'',\beta'} \bg'_{m-p} W_{p+\alpha'',\beta'}
                                  \overline W_{\alpha'',\beta'-q} \bg''_{q+n}
 \; = \; \sum_{p,q} \bg'_{m-p} \omega_{p,-q} \bg''_{q+n} .
\end{eqnarray*}
This implies (\ref{g'nug''}). Equation (\ref{g'g''}) follows from (\ref{||mu=omega00}). \qed

We put
$$
  c_k = \max_{j\in\mZ_J} |W_{k+j,j}|, \quad
  \|W\|_{\DT} = \sum_{k\in\mZ_J} c_k .
$$
Then
\begin{equation}
\label{DT_finite}
  \sum_{m\in\mZ_J} |\omega_{m,n}| \le \|W\|_\DT^2, \quad
  \sum_{n\in\mZ_J} |\omega_{m,n}| \le \|W\|_\DT^2 .
\end{equation}

\begin{dfn}
\label{dfn:f>0}
We say that $f=\sum_{j\in\mZ_J} f_j \eta^{jx}$ is nonnegative ($f\ge 0$) if the numbers $f(x)$ are real and nonnegative for all $x\in\mZ_J$.
\end{dfn}

\begin{lem}
\label{lem:f>0}
If $\calX = \mZ_J$ then the following three statements are equivalent.

(1) $f\ge 0$,

(2) $f = g\overline g$ for some $g\in\calH$,

(3) $f_k = \sum_{s\in\mZ_J} g_{k+s} \overline g_s$.
\end{lem}

{\it Proof}. Equivalence of statements (2) and (3) is obvious. Equivalence of (1) and (2) follows from two simple facts.

(a) For any $g',g''\in\calH$ and any $x\in\mZ_J$\;
$(g'g'')(x) = g'(x) g''(x)$.

(b) The numbers $f(x)$, $x\in\mZ_J$ determine uniquely $f\in\calH$. \qed

\begin{lem}
\label{lem:omega=measure}
Suppose $f\ge0$ and
$\displaystyle F(x) = \sum_{m,n\in\mZ_J} \omega_{m,-n} f_n \eta^{mx}$. Then $F\ge 0$.
\end{lem}

{\it Proof}. By Lemma \ref{lem:f>0} $f=g\overline g$ for some $g\in\calH$. Then
$$
    F(x)
  = \frac1J \sum_{j,l,n,s\in\mZ_J}
         W_{l+m,j}\overline W_{l,j-n} g_{n+s} \overline g_s \eta^{mx}.
$$
By using the change of summation indices $(n,s)\mapsto (k,q)$, $k=j-n$, $q=s-k$, we obtain:
$$
    F(x)
  = \frac1J \sum_{j,l,k,q\in\mZ_J} W_{l+m,j} g_{j+q} \overline W_{l,k} \overline g_{k+q} \eta^{mx}
  = \sum_{q\in\mZ} F^{(q)}(x) \overline F^{(q)}(x),
$$
where
$$
    F^{(q)}(x)
  = \sum_{j,l\in\mZ_J} W_{l,j} g_{j+q} \eta^{lx} .
$$
Hence, $F\ge 0$.  \qed

Now suppose $W$ is a unitary operator:
\begin{equation}
\label{Wunitary}
  \frac1J \sum_{l\in\mZ_J} W_{l,j} \overline W_{l,k} = \delta_{j,k}, \qquad
  j,k\in\mZ_J.
\end{equation}
We put $\one_\calX = \sqrt{J} \in\calH$. Hence,
$$
    \one_{\calX}(x)
  = \sum_{k\in\mZ_J} (\one_\calX)_k \eta^{kx}, \qquad
  (\one_\calX)_k = \sqrt{J}\delta_{k,0}, \quad
  \langle\one_\calX , \one_\calX\rangle = 1.
$$

\begin{lem}
\label{lem:one_finite}
Suppose $W$ is unitary. Then
\begin{equation}
\label{1finite}
  \sum_{n\in\mZ_J} \omega_{m,-n} (\one_\calX)_n = (\one_\calX)_m, \quad
  \sum_{m\in\mZ_J} (\one_\calX)_m \omega_{-m,n} = (\one_\calX)_n.
\end{equation}
\end{lem}

{\it Proof}. By (\ref{nu_finite}) and (\ref{Wunitary})
$$
    \sum_{n\in\mZ_J} \omega_{m,-n} (\one_\calX)_n
  = \sum_{j,l,n\in\mZ_J} \frac{W_{l+m,j} \overline W_{l,j-n}}{J} \sqrt{J} \delta_{n,0}
  = \sum_{j,l\in\mZ_J} \frac{W_{l+m,j} \overline W_{l,j}}{\sqrt{J}}
  = (\one_\calX)_m
$$
The second equation (\ref{1finite}) can be obtained analogously.
\qed

\subsection{Measure associated with a regular operator}
\label{ssec:meas_reg}

We start from a simple remark. Take any $f\in L^2(\mT)$. Then $f\in L^1(\mT)$ and we may associate with $f$ the Radon measure $f d\mu$ on $\mT$ i.e., a linear functional on $C(\mT)$:
$$
  C(\mT) \ni \ph \mapsto \int_\mT  f\ph\, d\mu .
$$
We use this observation in the following lemma.

\begin{lem}
\label{lem:regular}
Suppose $W\in\calR(\mT)$ and for any small $\eps>0$ the function $f_\eps = f_\eps(x)$ satisfies
\begin{equation}
\label{feps}
     f_\eps = \pi_Y f_\eps \in L^2(\mT), \quad
     \|f_\eps\|^2 = 1, \quad
     Y = (a-\eps,a+\eps).
\end{equation}
Let $L_a$ be determined by (\ref{L_a}). Then there exists the weak limit\footnote
{i.e., we regard here $|L_a * f_\eps|$ as functionals on $C(\mT)$}
\begin{equation}
\label{nu}
     \lim_{\eps\searrow 0} |L_a * f_\eps|^2
  =  \nu, \qquad
     \nu = \nu(x,a)
  =  \sum_{m\in\mZ} v_m(a) e^{im(x-a)} .
\end{equation}
The limit (\ref{nu}) is independent of the choice of the family $f_\eps$, satisfying (\ref{feps}). For any $a\in\mT$
\begin{equation}
\label{dmu=nudx}
  d\widetilde\mu_{W,a} = \frac1{2\pi}\nu(x,a)\, dx
\end{equation}
is a (non-negative) measure on $\mT$ and its norm as a functional on $C(\mT)$ satisfies
\begin{equation}
\label{|mu|}
     \|\widetilde\mu_{W,a}\|_{C^*}
  =  \int_\mT d\widetilde\mu_{W,a}
 \le \bc^2
  =  \|W\|^2_{\DT}.
\end{equation}
\end{lem}

\begin{rem}
\label{rem:nu=sum}
By using (\ref{limv}) we obtain the equation
\begin{equation}
\label{nu=sum}
  \nu(x,a) = \sum_{m,n\in\mZ} \omega_{m,n} e^{imx+ina}.
\end{equation}
\end{rem}

\begin{cor}
Suppose $W\in\calR(\mT)$. Then by (\ref{dim=omega}) we have $\|W\|_\mu^2=\bT(W)=\omega_{0,0}$. Therefore by (\ref{nu=sum})
$$
  \|W\|_\mu^2 = \int_{\mT^2} \nu(x,a) \, da dx.
$$
\end{cor}

{\it Proof of Lemma \ref{lem:regular}}. For any $f_\eps$ satisfying (\ref{feps})
\begin{equation}
\label{e|Lf|}
      \Big| \frac1{2\pi} \int_\mT e^{-ikx} |L_a * f_\eps|^2\, dx \Big|
  \le \|\Conv_{L_a}\|^2 \le \bc^2, \qquad
      k\in\mZ.
\end{equation}
Hence absolute values of Fourier coefficients of any function $\nu_\eps(a,\cdot) = |L_a * f_\eps|^2$ do not exceed $\bc^2$. Existence of the weak limit (\ref{nu}) in the class of distributions with bounded Fourier coefficients is equivalent to the existence of limits for all Fourier coefficients:
\begin{equation}
\label{limFourier}
     \xi_m
   = \lim_{\eps\searrow 0} \xi_{\eps,m}, \qquad
     \xi_{\eps,m}
   = \frac1{2\pi} \int_\mT e^{-imx} |L_a * f_\eps|^2\, dx ,
\end{equation}
independence of these limits of $f_\eps$, and the equation $\xi_m = e^{-ima} v_m(a)$. The convergence (\ref{limFourier}) does not need to be uniform in $m$.

We put $f_\eps = \sum_k f_{\eps,k} e^{ikx}$. Below for brevity we skip the subscript $\eps$ and write $f_k$ instead of $f_{\eps,k}$. Then
$L_a * f_\eps = \sum_k w_k(a) f_k e^{ikx}$ and
$$
    \xi_{\eps,m}
  = \sum_{k\in\mZ} w_{k+m}(a) \overline w_k(a) f_{k+m} \overline f_k .
$$
By (\ref{e|Lf|}) $|\xi_{\eps,m}|\le\bc^2$ for any $m\in\mZ$.

Now we show that for any $m\in\mZ$ convergence (\ref{limFourier}) takes place. We put
\begin{equation}
\label{xi}
      \widetilde\xi_{\eps,m}
  =   \sum_{k\in\mZ} \sum_{|l'|,|l''|\le B}
        w_{k+m}(a) \overline w_k(a) \frac{f_{k+m+l'} \overline f_{k+l''}}{(2B+1)^2} e^{i(l'-l'')a}.
\end{equation}
Then by (\ref{|w(a)|<}) and (\ref{f-fe})
\begin{eqnarray}
\nonumber
       |\xi_{\eps,m} - \widetilde\xi_{\eps,m}|
 &\le& \sum_{k\in\mZ} \sum_{|l'|,|l''|\le B} |w_{k+m}(a)\overline w_k(a)|
           \frac{|f_{k+m}\overline f_k
                      - f_{k+m+l'}\overline f_{k+l''} e^{i(l'-l'')a}|}
                {(2B+1)^2} \\
\nonumber
 &\le& \sum_{k\in\mZ} \sum_{|l'|,|l''|\le B} \!\! \bc^2
           \frac{|f_{k+m} - f_{k+m+l'} e^{il'a}|\cdot |f_k|
                      + |f_{k+m+l'}|\cdot |f_{k} - f_{k+l''} e^{il'' a}|}
                {(2B+1)^2} \\
\nonumber
 &\le& \frac{2\bc^2}{2B+1}
          \sum_{|l|\le B} \|f_\eps\|\cdot \|f_\eps - e^{-il(x-a)} f_\eps\| \\
\label{xi-xi}
 &\le& \frac{2\bc^2}{2B+1} \|f_\eps\|^2 \sum_{|l|\le B} \frac{\eps^{3/2}}{\sqrt\pi} |l|
\;\le\: 2\eps^{3/2} B\bc^2 \|f_\eps\|^2 .
\end{eqnarray}
By using in (\ref{xi}) the notation $l=l'$, $n=k+l'$, and $s=l'-l''$, we obtain:
$$
     \widetilde\xi_{\eps,m}
  =  \sum_{n\in\mZ} \sum_{|l|,|s-l|\le B}
        w_{n+m-l}(a) \overline w_{n-l}(a) \frac{f_{n+m} \overline f_{n-s}}{(2B+1)^2} e^{isa} .
$$

If we fix $n\in\mZ$ and $s\in [-2B,2B]\cap\mZ$ in the last sum then $n-l\in I_{n,s}$, where
$$
    I_{n,s}
  = \left\{ \begin{array}{cc}
              [n-B,n+B-s] \cap\mZ &  \mbox{ if } s\ge 0, \\ {}
              [n-B-s,n+B] \cap\mZ &  \mbox{ if } s < 0.
            \end{array}
    \right.
$$
We also put
$$
      b_m
  = \frac{2B + 1 - |m|}{(2B+1)^2}, \qquad  |m| \le 2B+1.
$$
Then $\sum_{|s|\le 2B} b_s = 1$.  By (\ref{vIm})
$$
    \widetilde\xi_{\eps,m}
  = \sum_{k\in\mZ,\,|s|\le 2B} v_{I_{k,s},m}(a)\,
                  b_s f_{k+m} \overline f_{k-s} e^{isa} .
$$

Since $W$ is regular, then by Lemma \ref{lem:lim_v} for any $m\in\mZ$ and any $\sigma>0$ there exists $Q_m = Q_m(\sigma)$ such that
\begin{equation}
\label{v-v}
  \# I \ge Q_m \quad\mbox{implies}\quad
  |v_{I,m}(a) - v_m(a)| < \sigma.
\end{equation}

Then $\widetilde\xi_{\eps,m} = \Sigma_1 + \Sigma_2 + \Sigma_3$, where
\begin{eqnarray*}
      \Sigma_1
  &=& \sum_{k\in\mZ,\,|s|\le 2B} v_m(a)\, b_s f_{k+m} \overline f_{k-s} e^{isa}, \\
      \Sigma_2
  &=& \sum_{k\in\mZ,\,2B+1-|s|\ge Q_m} (v_{I_{k,s},m}(a) - v_m(a))\,
                                 b_s f_{k+m} \overline f_{k-s} e^{isa}, \\
      \Sigma_3
  &=& \sum_{k\in\mZ,\,2B+1-|s| < Q_m} (v_{I_{k,s},m}(a) - v_m(a))\,
                                 b_s f_{k+m} \overline f_{k-s} e^{isa} .
\end{eqnarray*}

By using Lemma \ref{lem:YfJ}, we obtain:
\begin{eqnarray}
\nonumber
     \Big|\Sigma_1 - e^{-ima} v_m(a) \|f_\eps\|^2 \Big|
 &=& \bigg| \sum_{|s|\le 2B} e^{-ima} v_m(a) b_s \sum_{k\in\mZ}
                    \big( f_{k+m}\overline f_{k-s} e^{i(m+s)a} - |f_k|^2 \big)
     \bigg| \\
\label{Sigma1}
&\le& 3\bc^2 (2B + |m|) \eps \|f_\eps\|^2 .
\end{eqnarray}

If $2B + 1 -|s| \ge Q_m$ then by using (\ref{v-v}) we can estimate $\Sigma_2$:
\begin{equation}
\label{Sigma2}
      |\Sigma_2|
  \le \sum_{k\in\mZ,\,2B+1-|s|\ge Q_m} \sigma b_s  f_{k+m}\overline f_{k-s}
  \le \sigma \|f\|^2 .
\end{equation}

We have:
\begin{eqnarray}
\nonumber
        |\Sigma_3|
 &\le&  2\bc^2 \sum_{k\in\mZ,\,2B+1-|s| < Q_m} b_s |f_{k+m}\overline f_{k-s}| \\
\label{Sigma3}
 &\le&  2\bc^2 \sum_{2B+1-|s| < Q_m} b_s \|f_\eps\|^2
\;\le\; 2\bc^2 \frac{(1 + Q_m)Q_m}{(2B+1)^2} \|f_\eps\|^2.
\end{eqnarray}

Combining estimates (\ref{xi-xi}), (\ref{Sigma1}), (\ref{Sigma2}), and (\ref{Sigma3}), we obtain:
$$
       \big| \xi_{\eps,m} - e^{-ima} v_m(a) \|f_\eps\|^2 \big|
  \le  \bigg( (2\eps^{3/2} B + 6\eps B + 3\eps |m|)\bc^2 + \sigma + 2\bc^2 \frac{(1 + Q_m)Q_m}{(2B+1)^2}
       \bigg) \|f_\eps\|^2.
$$
This implies existence of the limits (\ref{limFourier}). Since the functions $\nu_\eps(a,\cdot)$ are nonnegative, $d\widetilde\mu_{W,a}$ is a measure on $\mT$ for any $a\in\mT$. \qed

\subsection{The space $\ACF(\mT)$}

\begin{dfn}
We say that $f\in L^\infty(\mT)$ lies in the space $\ACF(\mT)$ if it has absolutely converging Fourier series.
\end{dfn}

For any $f\in\ACF(\mT)$ we put
$$
  \|f\|_{\DT} = \sum_{k\in\mZ} |f_k| .
$$
Then for any $f\in\ACF(\mT)$ we have: $\|f\|_{\DT} = \|\widehat f\|_{\DT}$ and the estimate (\ref{||f||<||f||_DT}).
Assertion {\bf (14)} implies the following

\begin{cor}
\label{cor:ACF}
The space $\big(\ACF(\mT),\|\cdot\|_\DT\big)$ is a commutative $\mC^*$-subalgebra in the $\mC^*$-algebra
$\big(\DT(\mT),\|\cdot\|_\DT\big)$.
\end{cor}

Suppose $\nu$ satisfies (\ref{nu=sum}) for some $W\in\calR(\mT)$. Consider the functions
\begin{equation}
\label{phipsi}
   \ph = \frac1{2\pi}\int_{\mT} \nu(x,\cdot)\, dx
   \quad{and}\quad
   \psi = \frac1{2\pi}\int_{\mT} \nu(\cdot,a)\, da .
\end{equation}

\begin{lem}
\label{lem:intnuACF}
Suppose $W\in\calR(\mT)$. Then $\ph,\psi\in\ACF(\mT)$. Moreover,
$$
  \|\ph\|_\DT \le \|W\|^2_\DT, \quad
  \|\psi\|_\DT \le \|W\|^2_\DT.
$$
\end{lem}

{\it Proof}. By (\ref{nu=sum})
$$
    \frac1{2\pi} \int_\mT \nu(x,a) \, dx
  = \sum_{n\in\mZ} \omega_{0,n} e^{ina} .
$$
By (\ref{sumomega}) $\sum |\omega_{0,n}| \le \|W\|_\DT^2$.

The case of the function $\psi$ is analogous. \qed

\subsection{The operator $\widehat g_1 W \widehat g_2$}

\begin{prop}
\label{prop:WgWWg_reg}
Suppose $W\in\calR(\mT)$, $g_2\in L^\infty(\mT)$, and $g_1, |g_2|^2\in\ACF(\mT)$. Then the operator $\widetilde W = \widehat g_1 W \widehat g_2$ is also regular and the corresponding coefficients $\widetilde\omega_{m,n}$ satisfy
\begin{equation}
\label{tildeomega_e_WgWWg}
    \sum_{m,n\in\mZ} \widetilde\omega_{m,n} e^{imx+ina}
  = |g_1(x)|^2 \nu(x,a) |g_2(a)|^2.
\end{equation}
\end{prop}

\begin{cor}
\label{cor:WgWWg}
Suppose $W\in\calR(\mT)$, $g_2\in L^\infty(\mT)$, and $g_1, |g_2|^2\in\ACF(\mT)$. Then by (\ref{omega=T}) and (\ref{dim=omega})
\begin{equation}
\label{gWg}
     \| \widehat g_1 W \widehat g_2 \|_\mu^2
   = \frac{1}{(2\pi)^2} \int_{\mT^2} |g_1(x)|^2 \nu(x,a) |g_2(a)|^2\, da dx .
\end{equation}
\end{cor}

Proposition \ref{prop:WgWWg_reg} follows from Lemmas \ref{lem:WWg_reg} and \ref{lem:dimWW} while Corollary \ref{cor:WgWWg} is a combination of Corollaries \ref{cor:dimWWg} and \ref{cor:dimWgW}.
\smallskip

We associate with any $W\in\calR(\mT)$ the measure $\mu_W$ on $\calX\times\calX$:
$$
  d\mu_W = \nu(x,a)\, \frac{dx da}{4\pi^2} .
$$

\section{A bistochastic operator generated by $\mu_W$}
\label{sec:bistochastic}

\begin{lem}
\label{lem:bistochastic}
Suppose $W$ satisfies (at least) one of the following conditions:

(1) $W = U_F$, where $F\in\Aut(\calX,\mu)$,

(2) $W$ is an operator on $\calH = L^2(\calX)$, $\calX = \mZ_J$,

(3) $W\in\calR(\mT)$.

Then the corresponding measure $\mu_W$ determines a bounded operator
$$
  \calW : L^1(\calX,\mu)\to L^1(\calX,\mu) , \qquad
  f \mapsto \calW f = \int_\calX \nu(\cdot,a) f(a)\, d\mu(a) .
$$
This operator satisfies the estimate
\begin{equation}
\label{norm_calW}
  \|\calW\|_{L^1\to L^1} \le \|W\|_\DT^2
\end{equation}
and moreover, has the following properties

(a). $\calW$ is nonnegative: $\calW f\ge 0$ whenever $0\le f\in L^1(\calX,\mu)$.

If $W$ is unitary\footnote{In case (1) this condition automatically holds.} then two more statements hold.

(b). $\calW\one_\calX = \one_\calX$.

(c). $\int_\calX \calW f(x) \, d\mu(x) = \int_\calX f(a) \, d\mu(a)$ for any $f\in L^1(\calX,\mu)$.
\end{lem}

Conditions (a)--(c) mean that $\calW$ is a bistochastic (doubly stochastic) operator.
\medskip

{\it Proof of Lemma \ref{lem:bistochastic}}. (1) In this case $d\mu_W(x,a) = \delta(x,F(a))\, d\mu(x) d\mu(a)$ (see (\ref{muUF})). Then by (\ref{int_delta})--(\ref{delta=delta}) for any $f\in L^1(\calX)$
\begin{equation}
\label{calWf}
    \calW f(x)
  = \int_\calX \delta(x,F(a)) f(a) \, d\mu(a)
  = \int_\calX \delta(F^{-1}(x),a) f(a) \, d\mu(a)
  = f\circ F^{-1}(x).
\end{equation}
Hence, $\|\calW f\|_1 = \|f\|_1$. Therefore $\calW$ is an isometry and condition (c) holds. Conditions (a) and (b) also follow from (\ref{calWf}).
\smallskip

(2) If $\calX = \mZ_J$ then by (\ref{nu_finite}) we have:
$$
    f(a)
  = \sum_{j\in\mZ} f_j \eta^{ja} , \quad
    \calW f(x)
  = \sum_{m,n\in\mZ} \omega_{m,-n} f_n \eta^{mx}.
$$
Hence, by (\ref{DT_finite})
$$
      \|\calW f\|_1
  \le \sum_{m,n\in\mZ_J} |\omega_{m,-n} f_n|
  \le \|W\|_\DT^2 \|f\|_1 .
$$

Condition (a) follows from Lemma \ref{lem:omega=measure}, while Conditions (b) and (c) from Lemma \ref{lem:one_finite}.
\smallskip

(3) In this case for any $f\in L^1(\mT)$
$$
      \|\calW f\|_1
  =  \frac1{(2\pi)^2} \int_\mT dx \int_\mT |\nu(x,a) f(a)|\, da
  =  \frac1{2\pi} \int_\mT \ph(a) |f(a)| \, da ,
$$
where $\ph$ is determined by (\ref{phipsi}). By Lemma \ref{lem:intnuACF}
$\|\ph\|_\infty\le \|\ph\|_\DT\le \|W\|_\DT^2$. This implies (\ref{norm_calW}).

To prove (a), we fix $0 \le f\in L^1(\mT)$ and take any nonnegative continuous function $g:\mT\to\mR$. Then
$\frac1{2\pi} \int_\mT g(x) \calW f(x) \, dx$ is the result of the action of the measure $\nu$ on the non-negative function $g(x) f(a)$. This result is nonnegative. Hence, $\calW f\ge 0$.

Now suppose $W$ is unitary. Then by (\ref{gWg}) for any $0\not\equiv g\in\ACF(\mT)$
$$
     \frac1{2\pi} \int_\mT |g(x)|^2 \calW \one_\mT(x) \, dx
  =  \frac1{(2\pi)^2} \int_{\mT^2} |g(x)|^2 \nu(x,a) \, dx da
  = \|\widehat g W\|_\mu^2
  = \|\widehat g\|_\mu^2
$$
(the last equation follows from Corollary \ref{cor:UWU}). The equation
$$
    \|\widehat g\|_\mu^2
  = \|g\|^2 = \frac1{2\pi} \int_\mT |g(x)|^2 \, dx
$$
implies (b).

It is sufficient to check condition (c) only for $f\ge 0$. First, consider the case when $0\le f\in L^\infty(\mT)$. Consider $g\ge 0$ such that $f=g^2$. Then by Corollary \ref{cor:WgWWg}
$$
    \frac1{2\pi} \int_\mT \calW f(x) \, dx
  = \frac1{(2\pi)^2} \int_{\mT^2} \nu(x,a) |g(a)|^2 \, dx da
  = \|W\widehat g\|_\mu^2
  = \|\widehat g\|_\mu^2
  = \frac1{2\pi} \int_{\mT} f(a) \, da .
$$

If $0\le f\in L^1(\mT)$ is unbounded, by using a cut off, for any $\eps > 0$ we have: $f=f_1+f_2$,
$0\le f_1\in L^\infty(\mT)$, $\|f_2\|_1 < \eps$. Then by (\ref{norm_calW})
$$
     \|\calW f - \calW f_1\|_1
  <  \|W\|_\DT^2 \eps.
$$
This implies
$$
      \bigg| \frac1{2\pi} \int_\mT \calW f(x)\, dx - \frac1{2\pi} \int_\mT \calW f_1(x)\, dx \bigg|
  \le \|W\|_\DT^2 \eps .
$$
The estimate
$$
    \frac1{2\pi} \int_\mT \calW f_1(x)\, dx
  = \frac1{2\pi} \int_\mT f_1(x)\, dx
  = \frac1{2\pi} \int_\mT f(x)\, dx + \Delta, \qquad
    |\Delta| \le \eps
$$
finishes the proof. \qed

\section{Product with $\widehat g$}
\label{sec:product}

In this section we prove Proposition \ref{prop:WgWWg_reg} and Corollary \ref{cor:WgWWg}.

\subsection{The operator $W \widehat g$}

\begin{lem}
\label{lem:WWg_reg}
Suppose $W\in\calR(\mT)$, $g\in L^\infty(\mT)$, and $|g|^2\in\ACF(\mT)$. Then
$\widetilde W = W\widehat g$ is also regular and the corresponding coefficients $\widetilde\omega_{m,n}$ satisfy
\begin{equation}
\label{tildeomega_e}
    \sum_{m,n\in\mZ} \widetilde\omega_{m,n} e^{imx+ina}
  = \nu(x,a) |g(a)|^2.
\end{equation}
\end{lem}

\begin{cor}
\label{cor:dimWWg}
Suppose $W\in\calR(\mT)$ and $g\in\ACF(\mT)$. Then by (\ref{omega=T}) and (\ref{dim=omega})
\begin{equation}
\label{||Wg||mu}
     \| W\widehat g \|_\mu^2
   = \frac{1}{(2\pi)^2} \int_{\mT^2} \nu(x,a) |g(a)|^2 \, da dx .
\end{equation}
\end{cor}

{\it Proof of Lemma \ref{lem:WWg_reg}}. We denote by $\widetilde\omega_{I,m,n}$ the quantities (\ref{limomega}), corresponding to the operator
$\widetilde W = W\widehat g$. Then
$$
     \widetilde\omega_{I,m,n}
  =  \frac1{\# I} \sum_{k,s,j\in\mZ,\,l\in I}
         W_{l+m,k} g_{k-j} \overline W_{l,s} \overline g_{s-j-n}
  =  \frac1{\# I} \sum_{k,q\in\mZ,\,l\in I}
         W_{l+m,k}\overline W_{l,k-q+n} \bg_q,
$$
where $\bg_q = \sum_{p\in\mZ} g_p \overline g_{p-q}$.
Note that $\bg_q$ is the $q$-th Fourier coefficient of the function $|g|^2$. Indeed,
$$
    \bg_q
  = \langle g, g e^{iqx} \rangle
  = \frac1{2\pi} \int_\mT |g(x)|^2 e^{-iqx}\, dx.
$$

By assumption of the lemma the series $\sum \bg_q$ absolutely converges. We put
$$
      \sum_{q\in\mZ} |\bg_q|
   =  \bc_g.
$$
For any $\sigma>0$ there exists $Q(\sigma)\in\mN$ such that
$\sum_{|q|>Q(\sigma)} |\bg_q| < \sigma$. Given a small $\sigma$ we have: $\widetilde\omega_{I,m,n} = \Omega_1 + \Omega_2$, where
\begin{eqnarray*}
     \Omega_1
 &=& \frac1{\# I} \sum_{k\in\mZ,\,l\in I,\, |q|\le Q(\sigma)}
           W_{l+m,k} \overline W_{l,k-q+n} \bg_q, \\
     \Omega_2
 &=& \frac1{\# I} \sum_{k\in\mZ,\,l\in I,\, |q| > Q(\sigma)}
           W_{l+m,k} \overline W_{l,k-q+n} \bg_q .
\end{eqnarray*}

In the sum $\Omega_1$ the index $q$ contains only a finite number of values. Hence for all sufficiently big $\# I$ we have:
$$
    \Omega_1
  = \sum_{|q|\le Q(\sigma)} \omega_{m,n-q} \bg_q + \Delta, \qquad
    |\Delta|
 \le \sum_{|q|\le Q(\sigma)} \sigma |\bg_q|
 \le \sigma\bc_g .
$$

Now we estimate $\Omega_2$. Let $c_k$ be the majorating sequence for $W$ with
$\sum c_k = \bc$. Then
\begin{eqnarray*}
       |\Omega_2|
 &\le& \frac1{\# I} \sum_{k\in\mZ,\,l\in I\,|q|>Q(\sigma)}
             c_{l+m-k} c_{l-k+q-n} |\bg_q| \\
 &\le& \sum_{j\in\mZ,\,|q|>Q(\sigma)}
             c_{m+j} c_{q-n+j} |\bg_q|
   =   \sum_{|q| > Q(\sigma)} \widetilde c_{m+n-q} \sigma |\bg_q| ,
\end{eqnarray*}
where $\widetilde c_s = \sum_j c_{s+j} c_j$, \; $\sum_s \widetilde c_s = \bc^2$. Hence
$$
  |\Omega_2| \le \bc^2 \sum_{|q| > Q(\sigma)} |\bg_q| \le \bc^2 \sigma.
$$
By (\ref{sumomega})
$$
      \bigg| \sum_{|q| > Q(\sigma)} \omega_{m,n-q} \bg_q \bigg|
  \le \bc^2 \sum_{|q| > Q(\sigma)} |\bg_q|
  \le \bc^2 \sigma .
$$
Therefore
$$
      \Big| \widetilde\omega_{I,m,n} - \sum_{q\in\mZ} \omega_{m,n+q} \bg_q \Big|
  \le \sigma(\bc_g + 2\bc^2) .
$$

Finally note that
\begin{eqnarray*}
     \nu(x,a) |g(a)|^2
 &=& \sum_{q,p\in\mZ} g_p\overline g_{p-q} e^{iqa}
      \sum_{m,k\in\mZ} \omega_{m,k} e^{imx+ika} \\
 &=& \sum_{m,n,q,p\in\mZ} \omega_{m,n-q} g_p\overline g_{p-q} e^{imx+ina} .
\end{eqnarray*}

These computations imply that for all sufficiently big $\# I$ the quantity $\widetilde\omega_{I,m,n}$ differs arbitrarily small from
$$
     \sum_{p,q\in\mZ} \omega_{m,n+q} g_p\overline g_{p-q}
  =  \frac1{(2\pi)^2} \int_{\mT^2} e^{-imx-ina} \nu(x,a) |g(a)|^2 \, dx da .
$$
\qed





\subsection{The operator $\widehat g W$}

\begin{lem}
\label{lem:dimWW}
Suppose $W\in\calR(\mT)$ and $g\in\ACF(\mT)$. Then $\widetilde W = \widehat g W$ is also regular and the corresponding coefficients $\widetilde\omega_{m,n}$ satisfy
\begin{equation}
\label{breveomega)}
    \sum_{m,n\in\mZ}\widetilde\omega_{m,n} e^{imx+ina}
  = |g(x)|^2 \nu(x,a).
\end{equation}
\end{lem}

\begin{cor}
\label{cor:dimWgW}
Suppose $W\in\calR(\mT)$ and $g\in\ACF(\mT)$. Then by (\ref{omega=T}) and (\ref{dim=omega})
\begin{equation}
\label{dimWgW}
   \|\widehat g W\|_\mu^2
 = \frac1{(2\pi)^2} \int_{\mT^2} |g(x)|^2 \nu(x,a)\, dx da.
\end{equation}
\end{cor}

{\it Proof of Proposition \ref{lem:dimWW}}. We denote by $\widetilde\omega_{I,m,n}$ the quantities (\ref{limomega}), corresponding to the operator
$\widetilde W = \widehat g W$. Then
\begin{eqnarray*}
     \widetilde\omega_{I,m,n}
 &=& \frac1{\# I} \sum_{k,s,j\in\mZ,\, l\in I}
           g_{l+m-k} W_{k,j} \overline g_{l-s} \overline W_{s,j+n} \\
 &=& \frac1{\# I} \sum_{k,s,j\in\mZ,\, l\in I}
           g_{m-k} \overline g_{-s} W_{k+l,j} \overline W_{s+l,j+n} \\
 &=& \frac1{\# I} \sum_{k,s,j\in\mZ,\, l\in I+s}
           g_{m-k} \overline g_{-s} W_{k-s+l,j} \overline W_{l,j+n} \\
 &=& \sum_{k,s\in\mZ}  g_{m-k} \overline g_{-s} \omega_{I+s,k-s,n}.
\end{eqnarray*}
We define $\bc_g = \sum_{j\in\mZ} |g_j|$.
Given $\sigma>0$ we take $M = M(\sigma) > |m|$ such that
\begin{equation}
\label{sumg<sigma}
   \sum_{|j| \ge (M - |m|) / 2} |g_j| < \sigma.
\end{equation}
We take $N = N(\sigma,M)$ such that for any $j$, $|j|\le M$ and any interval $J\subset\mZ$, $\# J > N$
\begin{equation}
\label{om-om}
  |\omega_{J,j,n} - \omega_{j,n}| < \sigma.
\end{equation}
Then $\widetilde\omega_{I,m,n} = \Sigma_0 + \Sigma_1 + \Sigma_2 + \Sigma_3$, where
\begin{eqnarray*}
     \Sigma_0
 &=& \sum_{k,s\in\mZ} g_{m-k} \overline g_{-s} \omega_{k-s,n} , \\
     \Sigma_1
 &=& \sum_{|k-s|\ge M} g_{m-k} \overline g_{-s} \omega_{k-s,n} , \\
     \Sigma_2
 &=& - \sum_{|k-s| < M}
         g_{m-k} \overline g_{-s} (\omega_{I+s,k-s,n} - \omega_{k-s,n}) , \\
     \Sigma_3
 &=& \sum_{|k-s|\ge M} g_{m-k} \overline g_{-s} \omega_{I,k-s,n} .
\end{eqnarray*}

First, we transform the sum $\Sigma_0$:
\begin{eqnarray}
\nonumber
     \Sigma_0
 &=& \frac1{(2\pi)^2} \int_{\mT^2} e^{-imx-ina}
       \sum_{k,s\in\mZ} g_{m-k} \overline g_{-s}
            e^{i(m-k+s)x} \omega_{k-s,n} e^{i(k-s)x+ina}\, dxda \\
\label{Sig0}
 &=& \frac1{(2\pi)^2} \int_{\mT^2} e^{-imx-ina}
       |g(x)|^2 \nu(x,a) \, dxda .
\end{eqnarray}

We estimate $\Sigma_1$ by using (\ref{omegaI}) and (\ref{sumg<sigma}):
\begin{equation}
\label{Sig1}
       |\Sigma_1|
  \le  \bc^2 \sum_{|k-s| \ge M} |g_{m-k}| |g_{-s}|
  \le 2\bc^2 \sum_{\beta\in\mZ, |\alpha| \ge (M-|m|)/2} |g_\alpha| |g_\beta|
  \le  2\bc^2 \bc_g \sigma.
\end{equation}

We have by (\ref{om-om}):
\begin{equation}
\label{Sig2}
       |\Sigma_2|
  \le  \sigma \sum_{k,s\in\mZ} |g_{m-k}| |g_{-s}|
   =   \sigma\bc_g^2.
\end{equation}

The sum $\Sigma_3$ is estimated in the same way as $\Sigma_1$:
\begin{equation}
\label{Sig3}
       |\Sigma_3|
  \le  \bc^2 \sum_{|k-s| \ge M} |g_{m-k}| |g_{-s}|
  \le  2\bc^2\bc_g\sigma.
\end{equation}

Combining (\ref{Sig0}), (\ref{Sig1}), (\ref{Sig2}), and (\ref{Sig3}), we see that
$$
    \bigg| \widetilde\omega_{I,m,n}
        - \frac1{(2\pi)^2} \int_{\mT^2} e^{-imx-ina} |g(x)|^2 \nu(x,a) \, dxda
    \bigg|
  \le (4\bc^2\bc_g + \bc_g^2) \sigma.
$$
This implies (\ref{breveomega)}). \qed


\section{Weaker assumptions on $g$}

Unfortunately we do not know if equation (\ref{gWg}) remains valid for $W\in\calR(\mT)$ and arbitrary
$g_1,g_2\in L^\infty(\mT)$.
In this section we present two partial results in this direction. More precisely, we show that equations (\ref{||Wg||mu}) and (\ref{dimWgW}) hold for $g$ lying in sets which are larger than declared in Corollaries \ref{cor:dimWWg} and \ref{cor:dimWgW}, in particular, for $g$ equal to indicators of intervals.

\subsection{Computation of $\|W \widehat g\|_\mu$}

\begin{lem}
\label{lem:Whatg}
Suppose $W\in\calR(\mT)$ and $g\in L^\infty(\mT)$. Then equation (\ref{||Wg||mu}) remains valid.
\end{lem}

{\it Proof}. For any $\sigma > 0$ there exists $g_*\in\ACF(\mT)$ such that
$$
  \|g_*\|_\infty \le \|g\|_\infty \quad\mbox{and}\quad
  \|g - g_*\| < \sigma .
$$
Indeed, it is sufficient to take $\widetilde g$ equal to a finite, but sufficiently long part of the Fourier series of $g$, so that
$\big| \|g\|^2 - \|\widetilde g\|^2 \big| < \sigma$. Then we define $g_*$ as a cut off of $\widetilde g$:
$$
    g_*(x)
  = \left\{\begin{array}{ccc}
            \widetilde g(x) & \mbox{ if } & |\widetilde g(x)| \le \|g\|_\infty, \\
              \|g\|_\infty  & \mbox{ if } &  \widetilde g(x)  > \|g\|_\infty, \\
            - \|g\|_\infty  & \mbox{ if } &  \widetilde g(x)  < - \|g\|_\infty. \\
           \end{array}
    \right.
$$
The function $g_*(x)$ is Lipschitz. Hence, it lies in $\ACF(\mT)$.

By using the triangle inequality, (\ref{WW}) and Assertion {\bf (7)}, we have:
\begin{equation}
\label{g-g*1}
      \Big| \|W\widehat g\|_\mu - \|W\widehat g_*\|_\mu \Big|
  \le \|W(\widehat g - \widehat g_*)\|_\mu
  \le \|W\|\, \|(\widehat g - \widehat g_*)\|_\mu
  \le \sigma \|W\| .
\end{equation}

On the other hand, by (\ref{||Wg||mu})
\begin{eqnarray*}
\!\!\!\!
      \Big| \|W\widehat g_*\|_\mu^2
          - \frac1{(2\pi)^2} \int_{\mT^2} \nu(x,a) |g(a)|^2 \, dadx \Big|
 &=&  \frac1{2\pi} \Big| \int_{\mT} \ph(a) (|g_*(a)|^2 - |g(a)|^2) \, da \Big| ,    \\
      \ph
 &=&  \frac1{2\pi} \int_{\mT} \nu(x,\cdot) \, dx .
\end{eqnarray*}

By using the estimate $\|\ph\|_\infty \le \|\ph\|_\DT \le \|W\|_\DT^2$
(Lemma \ref{lem:intnuACF}), we obtain:
\begin{equation}
\label{g-g*2}
     \Big| \|W\widehat g_*\|_\mu^2
          - \frac1{(2\pi)^2} \int_{\mT^2} \nu(x,a) |g(a)|^2 \, dadx \Big|
 \le \|W\|_\DT^2 \sigma
\end{equation}
Hence, Lemma \ref{lem:Whatg} is a combination of (\ref{g-g*1}) and (\ref{g-g*2}). \qed

\subsection{Computation of $\|\widehat g W\|_\mu$}

\begin{dfn}
\label{dfn:AACF}
We say that $g\in L^\infty(\mT)$ is almost $\ACF(\mT)$ (the notation is $g\in\AACF(\mT)$) if there exists a continuous in
$\eps\in (-\eps_0,\eps_0)$ family of functions
$g_\eps\in L^\infty(\mT)$ such that the following conditions hold:

{\bf A1}. $g_0 = g$,

{\bf A2}. if $\eps\ne 0$ then $g_\eps\in\ACF(\mT)$,

{\bf A3}. $|g_\eps - g| \le |g_\eps - g_{-\eps}|$ if $0 \le \eps < \eps_0$,

{\bf A4}. $\lim_{\eps\to 0} \|g_\eps - g_{-\eps}\| = 0$.
\end{dfn}
A typical example of an $\AACF(\mT)$ function is the indicator of an interval (see Lemma \ref{lem:interval} below).

\begin{lem}
\label{lem:WfW}
Suppose $W\in\calR(\mT)$ and $g\in\AACF(\mT)$. Then equation (\ref{dimWgW}) remains valid.
\end{lem}

{\it Proof}. Let $g_\eps$ be the family from Definition \ref{dfn:AACF}.
By {\bf A3} and {\bf A1} for any $\eps\in (0,\eps_0)$ and any $X\in\calB$
$$
      0
  \le \|(\widehat g_\eps - \widehat g) W\pi_X\|^2
  \le \|(\widehat g_\eps - \widehat g_{-\eps}) W\pi_X\|^2 .
$$
This inequality implies
$$
      0
  \le \|(\widehat g_\eps - \widehat g) W\|_\mu^2
  \le \|(\widehat g_\eps - \widehat g_{-\eps}) W\|_\mu^2, \qquad
      \eps\in (0,\eps_0).
$$
By {\bf A2} and (\ref{dimWgW})
$\|(\widehat g_\eps - \widehat g) W\|_\mu^2 \le \Delta_\eps$,
$$
      \Delta_\eps
   =  \frac1{(2\pi)^2}\int_{\mT^2} |g_\eps(x) - g_{-\eps}(x)|^2 \nu(x,a) \, dxda
  \le \|g_\eps - g_{-\eps}\|^2 \|\psi\|_\infty , \quad
      \psi
   =  \frac1{2\pi} \int_\mT \nu(\cdot,a) \, da
$$
(see also (\ref{phipsi})). Note that $\|\psi\|_\infty \le \|\psi\|_\DT$. Hence by
{\bf A4} and Lemma \ref{lem:intnuACF} $\Delta_\eps \to 0$ as $\eps\to 0$.

By triangle inequality
$$
     \Big| \|\widehat g_\eps W\|_\mu^2
              - \frac1{(2\pi)^2} \int_{\mT^2} |g(x)|^2 \nu(x,a)\, dadx \Big|
  =  \big| \|\widehat g W\|_\mu^2 - \|\widehat g_\eps W\|_\mu^2
 \le \Delta_\eps .
$$
Hence equation (\ref{dimWgW}) follows. \qed


\subsection{Indicators of intervals}

\begin{lem}
\label{lem:interval}
For any interval $I\subset\mT$ the function $\one_I$ lies in $\AACF$.
\end{lem}

{\it Proof}. We put $g_0 = \one_I$. For $\eps\in (0,1)$ we define
\begin{eqnarray*}
      g_\eps(x)
  &=& \max \Big\{ 0 , 1 - \frac1{\eps} \dist(x,I) \Big\} , \\
      g_{-\eps}(x)
  &=& 1 - \max \Big\{ 0 , 1 + \frac1{\eps} \dist(x,\mT\setminus I) \Big\} .
\end{eqnarray*}

Conditions {\bf A1}--{\bf A4} from Definition \ref{dfn:AACF} also obviously hold.  Hence $\one_I\in\AACF$. \qed


\section{Coefficients $\omega_{m,n}$}

In this section $W\in\DT(\mT)$, $\{c_k\}_{k\in\mZ}$ is its majorating sequence, and $\bc = \sum c_k = \|W\|_\DT$.

\begin{lem}
\label{lem:omegabaromega}
For any $W\in\calR(\mT)$ and any $m,n\in\mZ$
$$
  \omega_{m,n} = \overline\omega_{-m,-n} .
$$
\end{lem}

{\it Proof}. This equation follows from the identity
$\omega_{I,m,n} = \overline\omega_{I+m,-m,-n}$ for any interval $I\subset\mZ$ .
\qed



We have the estimates
\begin{equation}
\label{omegaI}
\begin{array}{ccc}
\displaystyle
       \sum_{n\in\mZ} |\omega_{I,m,n}|
 &\le&
\displaystyle  \frac1{\# I} \sum_{j,n\in\mZ,\, l\in I} c_{l+m-j} c_{l-j-n}
   =   \bc^2 , \\ \displaystyle
       \sum_{m\in\mZ} |\omega_{I,m,n}|
 &\le&
\displaystyle  \frac1{\# I} \sum_{j,n\in\mZ,\, l\in I} c_{l+m-j} c_{l-j-n}
   =   \bc^2 .
\end{array}
\end{equation}
These estimates imply the following

\begin{cor}
\label{cor:omega}
If $W$ is regular then (compare with (\ref{DT_finite}))
\begin{equation}
\label{sumomega}
  \sum_{n\in\mZ} |\omega_{m,n}| \le \bc^2 , \quad
  \sum_{n\in\mZ} |\omega_{n,m}| \le \bc^2 \quad
  \mbox{ for any $m\in\mZ$}.
\end{equation}
\end{cor}

We put
$$
  P_M(m) = \{n\in\mZ : |n+m| \ge 2M \}.
$$

\begin{lem}
\label{lem:sumomega}
For any $m\in\mZ$
\begin{equation}
\label{sum_P}
  \sum_{n\in P_M(m)} |\omega_{m,n}| \le 2\bc \sum_{|k| \ge M} c_k.
\end{equation}
\end{lem}

{\it Proof}. For any integer interval $I$
\begin{eqnarray*}
       \sum_{n\in P_M(m)} |\omega_{I,m,n}|
 &\le& \frac1{\# I} \sum_{j\in\mZ,\, n\in P_M(m),\, l\in I} c_{l-j+m} c_{l-j-n} \\
  &=&   \sum_{k\in\mZ,\, n\in P_M(m)} c_{m+n-k} c_k
 \; =\; \sum_{k\in\mZ,\, s\in P_M(0)} c_{s-k} c_k \\
 &\le&  \sum_{k\in\mZ,\, |m|\ge M} c_m c_k
       + \sum_{m\in\mZ,\, |k|\ge M} c_m c_k
 \; = \; 2\bc \sum_{|k|\ge M} c_k.
\end{eqnarray*}
In the limit $\# I\to\infty$ we obtain (\ref{sum_P}).
\qed

\end{document}